\documentclass[10pt,a4paper]{article}
\usepackage{comment}
\usepackage{empheq}
\usepackage[T1]{fontenc}
\usepackage[mathscr]{euscript}
\usepackage{amsmath}
\usepackage{mathtools}
\usepackage{extarrows}
\usepackage{amsbsy}
\usepackage{amsfonts}
\usepackage{amssymb}
\usepackage{graphicx}
\usepackage{bbm}
\usepackage{listings}
\usepackage{color}
\usepackage[table]{xcolor}
\usepackage{xcolor}
\usepackage[left=1cm,right=2cm,top=2cm,bottom=2cm]{geometry}
\usepackage[latin1]{inputenc}
\usepackage{amsmath,amssymb,amsfonts}
\usepackage[english]{babel}
\usepackage{enumerate}
\usepackage{yhmath}
\usepackage{amsthm}
\usepackage{picture}
\usepackage{colortbl}
\usepackage{hyperref}
\usepackage{wrapfig} 
\usepackage{graphicx}
\usepackage{xfrac}  
\usepackage{faktor}  
\usepackage{pst-all} 
\usepackage{pstricks}
\usepackage{permute} 

\usepackage{array}
\usepackage{tikz}
\usepackage{comment} 
\usepackage{dsfont} 
\usepackage{amsmath}
\usepackage{amsfonts}
\usepackage{amssymb}

\definecolor{ashgrey}{rgb}{0.7, 0.75, 0.71}

\usetikzlibrary{arrows}
\usetikzlibrary{arrows,shapes,positioning}
\usetikzlibrary[patterns]
\usetikzlibrary{decorations.markings}
\usepackage{pgf,tikz}
\usepackage{pgfplots}
\tikzstyle directed=[postaction={decorate,decoration={markings,
    mark=at position .5 with {\arrow[]{stealth}}}}]
\usepackage{tikz-cd} 
\usetikzlibrary{cd}
\usetikzlibrary{calc}
\usetikzlibrary{trees,positioning,arrows,fit,shapes}
\usetikzlibrary{matrix} 
\usetikzlibrary{angles,quotes,calc}
\definecolor{gris-clair}{rgb}{0.8,0.8,0.8}

\tikzstyle directed=[postaction={decorate,decoration={markings, 
mark=at position .65 with {\arrow{latex}}}}]

\usepackage{tikz-cd}

\tikzset{
  symbol/.style={
    draw=none,
    every to/.append style={
      edge node={node [sloped, allow upside down, auto=false]{$#1$}}}
  }
}

\usepackage[all]{xy}
\def\dar[#1]{\ar@<1pt>[#1]\ar@<-1pt>[#1]}
\entrymodifiers={!!<0pt,0.7ex>+}

\usepackage{amsthm} 

\theoremstyle{definition}

\theoremstyle{definition}
\title{An easy way to find solutions of the Diophantine equation $A^{3}+B^{3}=C^{3}+D^{3}$}
\author{D. FOSSE, MSc. Physics\\
dominique.fosse@a3.epfl.ch}\date{\vspace{-5ex}}
\begin{document}
\maketitle
In \cite{Hirschhorn}, it is shown how the evenness of the function $f(x):=(x^{2}+16x-21)^{3}+(2 x^{2}-4x+42)^{3}$ provides a solution to the diophantine equation $A^{3}+B^{3}=C^{3}+D^{3}$. Let's generalize this by considering this time $g(x):=(a_{2}x^{2}+a_{1}x+a_{0})^{3}+(b_{2}x^{2}+b_{1}x+b_{0})^{3}$. Expanding and simplifying $g(x)-g(-x)$ gives:
\begin{equation}\label{eq-1}
g(x)-g(-x)=6(a_{1}a_{2}^{2}+b_{1}b_{2}^{2})x^{5}+2(6a_{0}a_{1}a_{2}+6b_{0}b_{1}b_{2}+a_{1}^{3}+b_{1}^{3})x^{3}+6(a_{0}^{2}a_{1}+b_{0}^{2}b_{1})x
\end{equation}
All we have to do is to annihilate the three coefficients of this last polynomial. We want to stay with rational expressions
, so we start by picking up $b_{2}$ from the coefficient of $x^{3}$. This coefficient is zero when $b_{2}=-\tfrac{6 a_{0}a_{1}a_{2}+a_{1}^{3}+b_{1}^{3}}{6 b_{0}b_{1}}$ and also, the coefficient of $x$ in \eqref{eq-1} is zero when $b_{1}=-\tfrac{a_{0}^{2}a_{1}}{b_{0}^{2}}$. Substitute this expression of $b_{1}$ in the previous expression of $b_{2}$ to obtain:
\begin{equation}\label{eq-2}
b_{1}=-\frac{a_{0}^{2}a_{1}}{b_{0}^{2}}\qquad\textnormal{and}\qquad
b_{2}=\frac{(6 a_{0}a_{2}+a_{1}^{2})b_{0}^{6}-a_{1}^{2}a_{0}^{6}}{6 b_{0}^{5}a_{0}^{2}}
\end{equation}
Let's substitute those expressions in the coefficient of $x^{5}$ in \eqref{eq-1}, that is, at this stage, the last one not to be null; we get $6(a_{1}a_{2}^{2}+b_{1}b_{2}^{2})=\tfrac{a_{1}^{3}(b_{0}^{6}-a_{0}^{6})(a_{0}^{6}a_{1}^{2}-12 a_{0}a_{2}b_{0}^{6}-a_{1}^{2}b_{0}^{6})}{6 a_{0}^{2}b_{0}^{12}}$. Of course, this term is null when $a_{0}=b_{0}$ or when $a_{1}=0$ but we also see that $a_{2}$ is of degree one on the other factor and that provides a rational expression of it: $a_{0}^{6}a_{1}^{2}-12 a_{0}a_{2}b_{0}^{6}-a_{1}^{2}b_{0}^{6}\implies a_{2}=\tfrac{a_{1}^{2}(a_{0}^{6}-b_{0}^{6})}{12 a_{0} b_{0}^{6}}$. We replace this in the expression of $b_{2}$ in \eqref{eq-2} and we get finally $b_{2}=-\tfrac{a_{1}^{2}(a_{0}^{6}-b_{0}^{6})}{12 a_{0}^{2} b_{0}^{5}}$. We can now conclude that
\begin{equation}\label{eq-3}
h(x):=\left(\left(\frac{a_{1}^{2}(a_{0}^{6}-b_{0}^{6})}{12 a_{0} b_{0}^{6}}\right)x^{2}+a_{1}x+a_{0}\right)^{3}+\left(-\left(\frac{a_{1}^{2}(a_{0}^{6}-b_{0}^{6})}{12 a_{0}^{2} b_{0}^{5}}\right)x^{2}-\left(\frac{a_{0}^{2}a_{1}}{b_{0}^{2}}\right)x+b_{0}\right)^{3}
\end{equation}
is an even function. Then $h(x)=h(-x)$ allows us to get rid of the denominators in \eqref{eq-3}. We also see that $a_{1}$ and $x$ behave exactly the same way; so we can drop one of those two variables, say $a_{1}$. Changing $x$ into $\tfrac{x}{y}$ and rearranging the terms provides a parametrization of $A^{3}+B^{3}=C^{3}+D^{3}$ in terms of binary quadratic forms; i.e.:
\begin{equation}\label{eq-4}
h_{2}(x):=\left(12 q^{3}p^{6}y^{2}+12 q^{2}p^{6}xy+q(q^{6}-p^{6})x^{2}\right)^{3}+\left(12 q^{2}p^{7}y^{2}-12 q^{4}p^{4}xy-p(q^{6}-p^{6})x^{2}\right)^{3}
\end{equation}
is an even function for the variable $x$ (after having renamed $q:=a_{0}$ and $p:=a_{1}$). The last step, for the sake of concision, is to change $y$ into $\tfrac{y}{2q p^{3}}$ to get finally:
\begin{empheq}[box=\fbox]{align}\label{eq-5}
&{\big(3q y^{2}+6q p^{3}xy-q(p^{6}-q^{6})x^{2}\big)}^{3}
+{\big(3p y^{2}-6p q^{3}xy+p(p^{6}-q^{6})x^{2}\big)}^{3}\nonumber\\
=&{\big(3q y^{2}-6q p^{3}xy-q(p^{6}-q^{6})x^{2}\big)}^{3}
+{\big(3p y^{2}+6p q^{3}xy+p(p^{6}-q^{6})x^{2}\big)}^{3}\nonumber
\end{empheq}

\end{document}